\documentclass[12pt]{article}
\usepackage{amssymb,amsmath,mathrsfs,amsthm,latexsym}
\usepackage{hyperref}
\usepackage{graphicx}
\usepackage{color}
\usepackage[OT2,OT1]{fontenc}
\usepackage{upquote}
\usepackage[numbers,sort&compress]{natbib}
\usepackage{alphabeta}
\usepackage{alltt}

\begin{document}

\title{\Large\bf Efficient application of the Voigt functions in the Fourier transform}

\author{
\normalsize\bf Sanjar M. Abrarov, Rehan Siddiqui, Rajinder K. Jagpal \\
\normalsize\bf and Brendan M. Quine}

\date{June 24, 2025}
\maketitle

\begin{abstract}
In this work, we develop a method for rational approximation of the Fourier transform (FT) based on the real and imaginary parts of the complex error function
\[
w(z) = e^{-z^2}(1 - {\rm{erf}}(-iz)) = K(x,y) + iL(x,y), \qquad z = x + iy,
\]
where $K(x,y)$ and $L(x,y)$ are known as the Voigt and imaginary Voigt functions, respectively. In contrast to our previous rational approximation of the FT, the expansion coefficients in this method are not dependent on the values of a sampled function. As the values of the Voigt functions remain the same, this approach can be used for rapid computation with help of look-up tables. Mathematica codes with some examples are presented.
\vspace{0.2cm}
\\
\end{abstract}

\section{Introduction}

The forward and inverse Fourier transforms (FTs) can be defined in a symmetric form as~\cite{Hansen2014, Bracewell2000}
\begin{equation}
\label{FT1} 
{\cal F}\left\{f(t)\right\}(\nu) = \int_{-\infty}^\infty {f(t)e^{-2\pi i\nu t}dt} = {\hat f}(\nu)
\end{equation}
and
\begin{equation}
\label{IFT} 
{\cal F}^{-1}\left\{\hat f(\nu)\right\}(t) = \int_{-\infty}^\infty {{\hat f}(\nu)e^{2\pi i\nu t}d\nu} = f(t),
\end{equation}
respectively. In~this work, we will consider only the forward FT \eqref{FT1} since due to symmetric form the approximations for the inverse FT (IFT) \eqref{IFT} can be readily obtained from the forward FT by change of the~variables.

Equation \eqref{IFT} shows how the function $\hat{f}(\nu)$ is transformed from $\nu$-domain ({\it e.g.} frequency domain) back to its original $t$-domain ({\it e.g.} time domain). More generally, however, the~IFT is also applicable for transformation from $t$-domain to $\nu$-domain. This can be written conveniently in a shorthand notation as
\[
{\cal F}^{-1}{\Omega(t)}(\nu) = {\hat \Omega}^{-1}(\nu),
\]
where $\Omega(t)$ is an arbitrary function with $t$-domain.

There are relations between a function $f(t)$ and its even $f_{even}(t)$ and odd $f_{odd}(t)$ components. In~particular, the~even and odd functions can be readily generated by using the following relations
\begin{equation}
\label{EFI} 
f_{even}(t) = \frac{f(t) + f(-t)}{2}
\end{equation}
and
\begin{equation}
\label{OFI} 
f_{odd}(t) = \frac{f(t) - f(-t)}{2}
\end{equation}
such that
\[
f(t) = f_{even}(t) + f_{odd}(t).
\]
Therefore, due to the linearity of the FT we can also state that
\begin{equation}
\label{FTI} 
{\hat f(\nu)} = {\hat f_{even}}(\nu) + {\hat f_{odd}}(\nu).
\end{equation}

Using Euler’s identity
\[
e^{ix} = \cos(x) + i\sin(x),
\]
we can express Equation \eqref{FT1} as
\begin{equation}
\label{FTTF} 
{\cal F}\left\{f(t)\right\}(\nu) = \int_{-\infty}^\infty {f(t)\left(\cos(2\pi \nu t) - i\sin(2\pi\nu t)\right)}\,dt.
\end{equation}

Since an even function satisfies the condition
\begin{equation}
\label{I4EF} 
f_{even}(t) = f_{even}(-t),
\end{equation}
it follows that
\[
\int_{-\infty}^\infty {f_{even}(t)\cos(2\pi\nu t)}\,dt = 2\int_0^\infty {f_{even}(t)\cos(2\pi\nu t)}\,dt
\]
and
\[
\int_{-\infty}^\infty {f_{even}(t)\sin(2\pi\nu t)}\,dt = 0.
\]

Since an odd function satisfies the condition
\begin{equation}
\label{I4OF} 
-f_{odd}(t) = f_{odd}(-t),
\end{equation}
we have
\[
\int_{-\infty}^\infty {f_{odd}(t)\cos(2\pi\nu t)}\,dt = 0
\]
and
\[
\int_{-\infty}^\infty {f_{odd}(t)\sin(2\pi\nu t)}\,dt = 2\int_0^\infty  {f_{odd}(t)\sin(2\pi\nu t)}\,dt.
\]

Consequently, from~these relations and Equations \eqref{FTI} and \eqref{FTTF} we get the following identities
\begin{equation}
\label{FT2} 
{\cal F}\left\{f_{even}(t)\right\}(\nu) = 2\int_0^\infty {f_{even}(t)\cos(2\pi\nu t)}\,dt = {\hat f}_{even}(\nu)
\end{equation}
and
\begin{equation}
\label{FT3} 
{\cal F}\left\{f_{odd}(t)\right\}(\nu) = -2i\int_0^\infty {f_{odd}(t)\sin(2\pi\nu t)}\,dt = {\hat f}_{odd}(\nu).
\end{equation}

Equations \eqref{FT2} and \eqref{FT3} are of primary importance since they will be used 
in derivation of the new approximation of the~FT.

In a recent publication, we developed a new methodology providing a rational approximation of the FT~\cite{Abrarov2021}. Specifically, the~FT of a function $f(t)$ can be approximated as a rational function consisting of low-order polynomials, $3 \times 4$, 
in~form
\begin{equation}
\label{RA} 
\begin{aligned}
{\hat f}(\nu) &\approx \sum_{m = 1}^M {\left(\frac{\alpha_m + \beta_m \nu^2}{\kappa_m + \lambda_m\nu^2 + \nu^4} - i\frac{\gamma_m\nu + \theta_m\nu^3}{\kappa_m + \lambda_m\nu^2 + \nu^4}\right)}
\\
&= \sum_{m = 1}^M {\frac{\alpha_m - i\gamma_m\nu + \beta_m\nu^2 - i\theta_m\nu^3}{\kappa_m + \lambda_m\nu ^2 + \nu^4}},
\end{aligned}
\end{equation}
where expansion coefficients are given by
\small
\[
\begin{aligned}
\alpha_m &= \frac{1}{8M\pi^4}\sum_{n = -N}^N {f_{even}(nh)e^{nh\sigma}\left(\mu_m^2 + \sigma^2\right)\left(\sigma\cos(nh\mu_m) + \mu_m\sin(nh\mu_m)\right)},
\\
\beta_m &= \frac{1}{2M\pi^2}\sum_{n = - N}^N {f_{even}(nh)e^{nh\sigma}\left(\sigma\cos(nh\mu_m) - \mu_m\sin(nh\mu_m)\right)},
\end{aligned}
\]
\[
\begin{aligned}
\gamma_m &= \frac{1}{4M\pi^3}\sum_{n = -N}^N {f_{odd}(nh)e^{nh\sigma}\left(\left(\sigma^2 - \mu_m^2\right)\cos(nh\mu_m) + 2\sigma\mu_m\sin(nh\mu_m)\right)},
\\
\theta_m &= \frac{1}{M\pi}\sum_{n = -N}^N {f_{odd}(nh)e^{nh\sigma}\cos(nh\mu_m)},
\\
\kappa_m &= \frac{1}{16\pi^4}\left(\mu_m^2 + \sigma^2\right)^2,
\\
\lambda_m &= \frac{1}{2\pi^2}\left(\sigma^2 - \mu_m^2\right),
\\
{\mu _m} &= \frac{\pi(m - 1/2)}{Mh}
\end{aligned}
\]
\normalsize
and $h$, $\sigma$ are fitting~parameters.

It should be noted that in our previous publication~\cite{Abrarov2015}, we {\it de facto} applied a rational approximation of the FT. Therefore, a~new method of rational approximation, described in our paper~\cite{Abrarov2021}, is just a generalization of a sampling and integration technique that we proposed earlier in~\cite{Abrarov2015} to derive rapid and high-accuracy rational approximations of the complex error function. Specifically, in~our earlier publication~\cite{Quine2013}, we found a very useful product-to-sum identity
\begin{equation}
\label{SFI} 
\prod_{n = 1}^k\cos\left(\frac{x}{2^n}\right) = \frac{1}{2^{k - 1}}\sum_{n = 1}^{2{^{k - 1}}}\cos\left(\frac{n - 1/2}{2^{k - 1}}x\right).
\end{equation}

Since according to Euler's formula for the sinc function (this formula is also attributed to Vi\`ete) \cite{Kac1959, Gearhart1990}
\[
{\rm sinc}(x) = \frac{\sin(x)}{x} = \prod_{n = 1}^\infty\cos\left(\frac{x}{2^n}\right),
\]
from Equation \eqref{SFI} it follows that~\cite{Abrarov2015}
\[
{\rm sinc}(x) = \lim_{K \to \infty}\frac{1}{K}\sum_{n = 1}^K\cos\left(\frac{n - 1/2}{K}x\right),
\]
where $K = 2^{k - 1}$. Consequently, from~this limit we can find the following approximation
\begin{equation}
\label{ICE} 
{\rm sinc}(x) \approx \frac{1}{K}\sum_{n = 1}^K\cos\left(\frac{n - 1/2}{K}x\right), \qquad\qquad K\gg 1
\end{equation}
that can be used efficiently as a sampling function~\cite{Abrarov2015}. Since in this approximation the integer $K$ is a finite number, we can regard this sampling function as an incomplete cosine expansion of the sinc~function.

Applying the sampling function \eqref{ICE} together with a new method of rational approximation based on integration with exponential decay multiplier that we developed in \cite{Abrarov2021,Abrarov2015}, we can obtain the rational approximation \eqref{RA} of the FT. Some examples of rational approximation \eqref{RA} of the FT are shown in our work~\cite{Abrarov2021} and can be visualized by running a provided Matlab code (Matlab code can be copy-pasted from this link: \url{https://arxiv.org/pdf/2001.07533}).

There are many methods for rational approximations such as Pad\'e approximation~\cite{Baker1961, Brezenski1996}, minimax approximation~\cite{Filip2018, Nakatsukasa2020}, Remez algorithm~\cite{Pachon2009, Hofreither2021}, and so on. However, to~the best of our knowledge, a~method of rational approximation that we developed in~\cite{Abrarov2021} is new and has never been reported in 
scientific~literature.

The small-order polynomials $3 \times 4$ in quotients of rational approximation \eqref{RA} may be advantageous for numerical analysis with any kind of computations including matrix manipulations, integrations and differentiations. However,~rational approximation \eqref{RA} requires re-computation of the four expansion coefficients $\alpha_m$, $\beta_m$, $\gamma_m$, and $\theta_m$ every time the sampled function $f(t)$ changes. In~this work, we propose a new method for rational approximation of the 
FT based on a sum of the real and imaginary parts of the complex error function~\cite{Abramowitz1972, Armstrong1972, Schreier1992}. Such a representation of rational approximation of the FT does not require re-computation when shape of the sampled function $f(t)$ changes.

\section{Preliminaries}

The complex error function, also commonly known as the Faddeeva function, can be defined as \cite{Schreier1992, Armstrong1967}

\begin{equation}
\label{CEF1} 
w(z) = e^{-{z^2}}\left(1 + \frac{2i}{\sqrt\pi}\int_0^z {e^{t^2}dt}\right),
\end{equation}
where $z = x + iy$. Comparing equation \eqref{CEF1} with definition of the error function \cite{Abramowitz1972}
\[
{\rm{erf}}(z) = \frac{2}{\sqrt\pi}\int_0^z {e^{-{t^2}}}dt
\]
one can see that the complex error function can be expressed as \cite{Abramowitz1972, Zaghloul2012}
\begin{equation}
\label{CEF2} 
w(z) = e^{-{z^2}}\left[1 - {\rm erf}(-iz)\right].
\end{equation}
Therefore, the complex error function $w(z)$ can be considered as a reformulation of the error function ${\rm erf}(z)$.

Complex error function satisfies the following relation
\begin{equation}
\label{I4CEF} 
w(-z) = 2e^{-{z^2}} - w(z).
\end{equation}

It can be shown that the complex error function is a solution of the following differential equation
\[
w'(z) + 2zw(z) = \frac{2i}{\sqrt\pi},
\]
with initial condition
\[
w(0) = 1.
\]

The complex error function $w(z)$ is closely related to the complex probability function \cite{Armstrong1972, Schreier1992}
\begin{equation}
\label{CPF} 
W(z) = \frac{i}{\pi}\int_{-\infty}^\infty {\frac{e^{-t^2}}{z - t}}dt.
\end{equation}

There is a direct relationship between these two functions. In particular, both functions are equal to each other on the upper half of the complex plane \cite{Armstrong1972, Schreier1992}
\[
w(z) = W(z), \qquad {\rm Im}[z] > 0.
\]

Separating the real and imaginary parts of the complex probability function \eqref{CPF} as
\[
W(z) = K(x,y) + iL(x,y),
\]
results in
\begin{equation}
\label{KF} 
K(x,y) = \frac{y}{\pi}\int_{-\infty}^\infty \frac{e^{-t^2}}{y^2 + (x - t)^2}dt
\end{equation}
and
\begin{equation}
\label{LF} 
L(x,y) = \frac{1}{\pi}\int_{-\infty}^\infty \frac{e^{-t^2}(x - t)}{y^2 + (x - t)^2}dt,
\end{equation}
respectively. 

The real part $K(x,y)$ of the complex probability function is known as the Voigt function that is widely used in Atmospheric Physics to describe absorption and emission of atmospheric molecules \cite{Berk2017, Pliutau2017, Pliutau2021}. 

The imaginary part $L(x,y)$ is also used in various fields of Physics and Engineering \cite{Balazs1969, Chan1986}. It does not have a specific name. However, following Zaghloul and Ali \cite{Zaghloul2012}, for the convenience we will also refer to the function $L(x,y)$ as the imaginary Voigt function.

It is not difficult to show that substituting the following identity \cite{Armstrong1972}
\begin{equation}
\label{I4I} 
\frac{y}{y^2 + (x - t)^2} = \int_0^\infty  e^{-yq}\cos\left[(x - t)q\right]\,dq, \qquad y > 0,
\end{equation}
into equation \eqref{KF}, we obtain \cite{Armstrong1972, Srivastava1987}
\[
\begin{aligned}
K(x,y) &= \frac{1}{\pi}\int_{-\infty}^\infty e^{-t^2}\frac{y}{y^2 + (x - t)^2}dt
\\
&= \frac{1}{\pi}\int_0^\infty {\int_{-\infty}^\infty e^{-t^2}e^{-yq}\cos\left[(x - t)q\right]dt\,dq}
\\
&= \frac{1}{\pi}\int_0^\infty e^{-yq}\int_{-\infty}^\infty e^{-t^2}\cos\left[(x - t)q\right]dt\,dq 
\end{aligned}
\]
and since
\[
\int_{-\infty}^\infty e^{-t^2}\cos\left[(x - t)q\right]\,dt = \sqrt\pi e^{-q^2/4}\cos(qx), 
\]
we can write
\begin{equation}
\label{VF} 
K(x,y) = \frac{1}{\sqrt\pi}\int_0^\infty e^{-t^2/4 - yt}\cos(xt)\,dt, \qquad y > 0.
\end{equation}

Similarly, substituting the identity 
$$
\frac{x - t}{y^2 + (x - t)^2} = \int_0^\infty e^{-yq}\sin(q(x - t))\,dq, \qquad y > 0,
$$
into equation \eqref{LF} we get \cite{Srivastava1987}
\[
\begin{aligned}
L(x,y) &= \frac{1}{\pi}\int_{-\infty}^\infty e^{-t^2}\frac{x-t}{y^2 + (x - t)^2}dt
\\
&= \frac{1}{\pi}\int_0^\infty {\int_{-\infty}^\infty e^{-t^2}e^{-yq}\sin\left[(x - t)q\right]dt\,dq}
\\
&= \frac{1}{\pi}\int_0^\infty e^{-yq}\int_{-\infty}^\infty e^{-t^2}\sin\left[(x - t)q\right]dt\,dq 
\end{aligned}
\]
and since
\[
\int_{-\infty}^\infty e^{-t^2}\sin\left[(x - t)q\right]\,dt = \sqrt\pi e^{-q^2/4}\sin(qx), 
\]
we have
\begin{equation}
\label{IVF} 
L(x,y) = \frac{1}{\sqrt\pi}\int_0^\infty e^{-{t^2}/4 - yt}\sin(xt)dt, \qquad y > 0.
\end{equation}

Sum of the equations above in terms of the real and imaginary parts yields \cite{Srivastava1987}
\[
\begin{aligned}
K(x,y) + iL(x,y) &= \frac{1}{\sqrt\pi}\int_0^\infty e^{-t^2/4 - yt}\left[\cos(xt) + i\sin(xt)\right]\,dt
\\
&= \frac{1}{\sqrt\pi}\int_0^\infty e^{-t^2/4 - yt}e^{ixt}\,dt
\\
&= \frac{1}{\sqrt\pi}\int_0^\infty e^{-t^2/4}e^{(-y - ix)t}\,dt
\\
&= e^{(y - ix)^2}\left[1 - {\rm erf}(y - ix)\right].
\end{aligned}
\]
We can see now that this equation is consistent with equation \eqref{CEF2} of the complex error function $w(z)$.

\section{Results and discussion}

\subsection{Methodology}

Consider the rectangular function (solitary rectangular function) that can be \mbox{defined as}
\[
f_r(t) = \left\{
\begin{aligned}
&1, &\quad  &-1/2 < t < 1/2,
\\
&1/2, &\quad &\left| t \right| = 1/2,
\\
&0, &\quad &{\rm{otherwise.}}
\end{aligned}
\right.
\]

This function can be expressed by the following limit:
\begin{equation}
\label{L4RF} 
f_r(t) = \lim\limits_{k \to \infty} \frac{1}{\left(2t\right)^{2k} + 1}.
\end{equation}
Consequently, we can approximate the rectangular function by taking sufficiently large value of the parameter $k$.

Figure~1 shows the rectangular function and its approximation at $k = 35$ by dashed red and light blue colors, respectively. As~we can see, at~$k = 35$ Equation \eqref{L4RF} approximates the rectangular function reasonably well. Therefore, we can use the following approximation:
\begin{equation}
\label{A4RF} 
f_r(t) \approx \frac{1}{(2t)^{70} + 1}
\end{equation}
and apply it to perform numerically the~FT.

\begin{figure}[ht]
\begin{center}
\includegraphics[width=18pc]{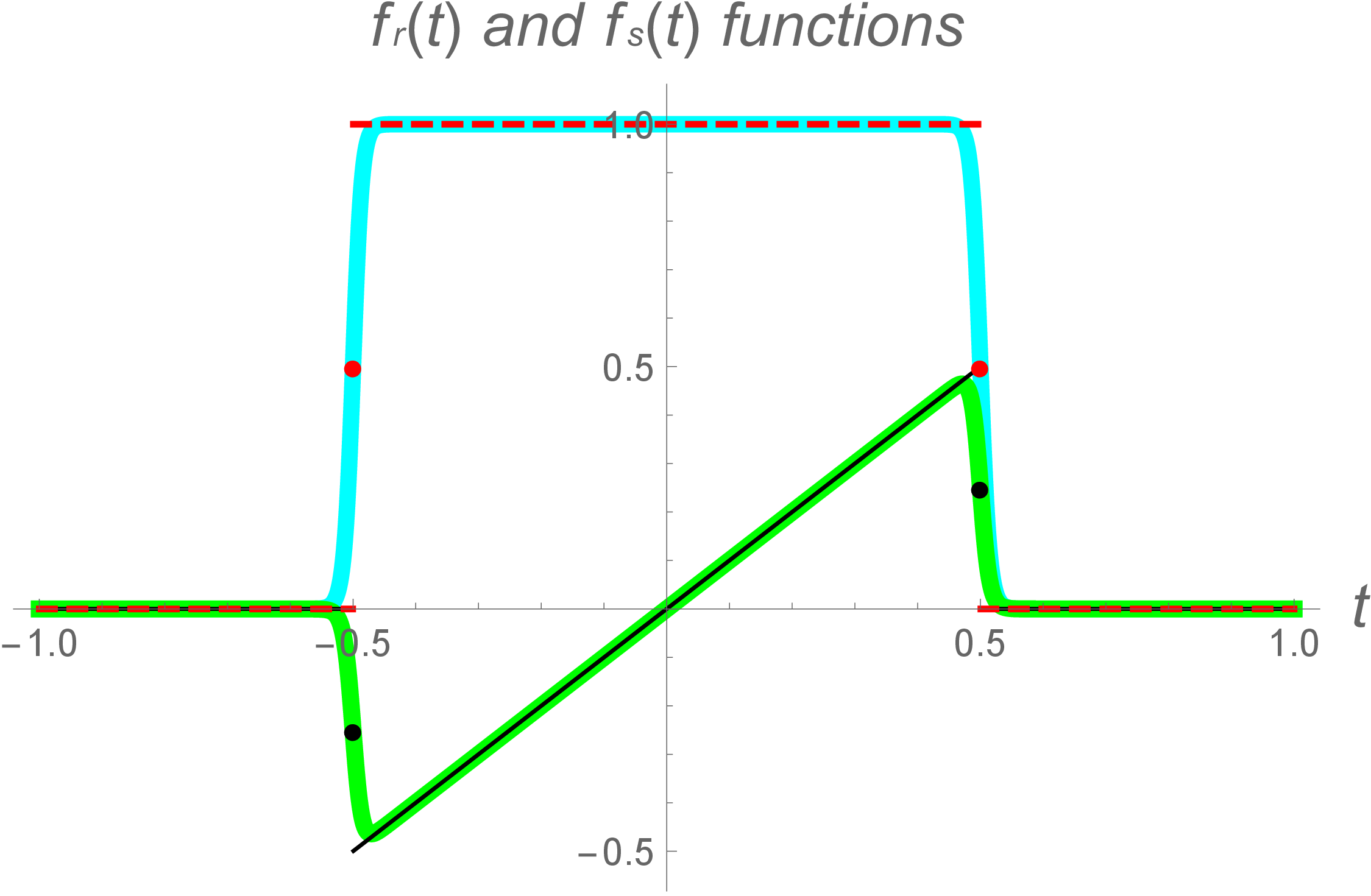}\hspace{2pc}%
\begin{minipage}[b]{28pc}
\vspace{0.3cm}
{\sffamily {\bf{Fig. 1.}} Rectangular and sawtooth functions and their approximations. Rectangular and sawtooth functions are shown by dashed red and solid black lines. Even $1/\left(\left(2t \right)^{70}+1\right)$ and odd $t/\left(\left(2t \right)^{70}+1\right)$ functions are shown by blue and green curves.}
\end{minipage}
\end{center}
\end{figure}

Previously, we have used the following sampling function:
\begin{equation}
\label{SF} 
s(t) = \frac{h}{c\sqrt \pi}{e^{-\left(\frac{t}{c}\right)^2}},
\end{equation}
where $h$ and $c$ are small fitting parameters, for~high-accuracy approximation of the complex error function~\cite{Abrarov2016}. Thus, applying this sampling function over the points $nh$ to approximation \eqref{A4RF}, the~rectangular function can be approximated as
\begin{equation}
\label{SA4RF} 
\begin{aligned}
f_r(t) &\approx \sum_{n = -N}^N {s(t - nh)f_r(nh)}
\\
&= \frac{h}{c\sqrt\pi}\sum_{n = -N}^N {e^{-(\frac{t - nh}{c})^2}f_r(nh)}
\\
&= \frac{h}{c\sqrt\pi}\sum_{n = -N}^N {e^{\left(2\frac{nh}{c^2}\right)t}e^{-\left(\frac{t}{c}\right)^2}e^{-\left(\frac{nh}{c}\right)^2}f_r(nh)}. 
\end{aligned}
\end{equation}

Figure~2 shows the approximations of the rectangular function at different fitting parameters $h$ and $c$. Specifically, the~green curve corresponds to $h = 0.065$, $c = 0.035$, and $N = 25$, the~red curve corresponds to $h = 0.05$, $c = 0.03$, and $N = 25$ while the light blue curve corresponds to $h = 0.02$, $c = 0.025$, and $N = 25$. The~rectangular function is also shown by black dashed~curve.

\begin{figure}[ht]
\begin{center}
\includegraphics[width=18pc]{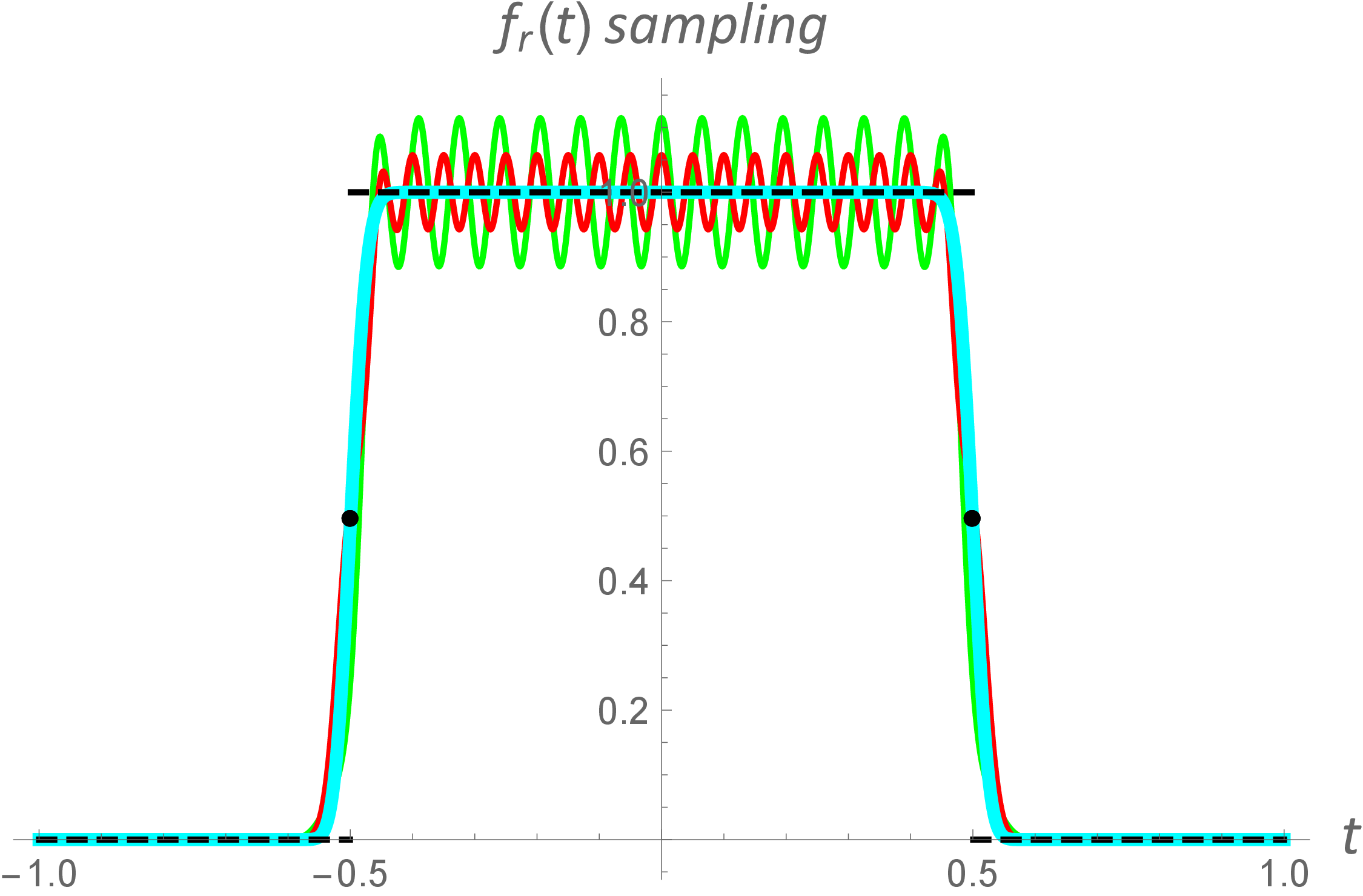}\hspace{2pc}%
\begin{minipage}[b]{28pc}
\vspace{0.3cm}
{\sffamily {\bf{Fig. 2.}} The rectangular function (dashed black curve) and its approximations by sampling at $h=0.065$, $c=0.035$, $N=25$ (green curve), $h=0.05$, $c=0.03$, $N=25$ (red curve) and $h=0.02$, $c=0.025$, $N=25$ (light blue curve).}
\end{minipage}
\end{center}
\end{figure}

Since the rectangular function $f_r(t)$ is even, we can use Equation \eqref{FT2} for the FT. Thus, the~substitution of approximation \eqref{SA4RF} into Equation \eqref{FT2} leads to
\[
\begin{aligned}
{\hat f}_r(\nu) &\approx 2\frac{h}{c\sqrt\pi}\int_{0}^\infty \sum_{n = -N}^N e^{-\left(\frac{t}{c}\right)^2 + \left(2\frac{nh}{c^2}\right)t}e^{-\left(\frac{nh}{c}\right)^2}f_r(nh)\cos(2\pi\nu t)dt
\\
&= 2\frac{h}{\sqrt\pi}\int_{0}^\infty \sum_{n = -N}^N e^{-t^2 + \frac{2nh}{c}t}e^{-\left(\frac{nh}{c}\right)^2}f_r(nh)\cos(2\pi\nu ct)dt
\\
&= \frac{h}{\sqrt\pi}\int\limits_{0}^\infty \sum_{n = -N}^N e^{-\frac{t^2}{4} - \frac{nh}{c}t}e^{ - \left(\frac{nh}{c} \right)^2}f_r(nh)\cos(\pi\nu ct)dt.
\end{aligned}
\]
Comparing the  approximation above with Equation \eqref{VF}, one can see that the FT of the rectangular function \eqref{L4RF} can be expressed in terms of the Voigt function
\begin{equation}
\label{FT4RF}
{\hat f}_r(\nu) \approx h\sum_{n = -N}^N e^{-\left(\frac{nh}{c}\right)^2}f_r(nh)K\left(\pi\nu c,\frac{nh}{c}\right).
\end{equation}

As a simple example for an odd function, we can consider the sawtooth function (solitary 
sawtooth function)
\[
f_s(t) = tf_r(t) = \left\{
\begin{aligned}
&t, &\quad  &-1/2 < t < 1/2,
\\
&-1/4, &\quad & t = -1/2,
\\
&1/4, &\quad & t = 1/2,
\\
&0, &\quad &{\rm{otherwise.}}
\end{aligned}
\right.
\]
Since this function can be express though the following limit
\begin{equation}
\label{L4SF} 
f_s(t) = t\lim_{k\to\infty} \frac{1}{(2t)^{2k} + 1},
\end{equation}
we can approximate it by multiplying $t$ with Equation \eqref{A4RF} as follows
\begin{equation}
\label{A4SF} 
f_s(t) \approx t\,\frac{1}{(2t)^{70} + 1}.
\end{equation}

Figure~1 shows the sawtooth function \eqref{L4SF} and its approximation \eqref{A4SF} by black and light green curves, respectively.

Applying the sampling function \eqref{SF} over the points $nh$ to approximation \eqref{A4SF}, we obtain
$$
\begin{aligned}
f_s(t) &\approx \sum_{n = -N}^N {s(t - nh)f_s(nh)}
\\
&= \frac{h}{c\sqrt\pi}\sum_{n = -N}^N e^{-(\frac{t - nh}{c})^2}f_s(nh)
\\
&= \frac{h}{c\sqrt\pi}\sum_{n = -N}^N e^{\left(2\frac{nh}{c^2}\right)t}e^{-\left(\frac{t}{c}\right)^2}e^{-\left(\frac{nh}{c}\right)^2}f_s(nh). 
\end{aligned}
$$
Consequently, the~FT of the sawtooth function \eqref{L4SF} can be approximated as
\[
\begin{aligned}
{\hat f}_s(\nu) &\approx -2i\frac{h}{c\sqrt\pi}\int_{0}^\infty\sum_{n = -N}^N e^{-\left(\frac{t}{c}\right)^2 + \left(2\frac{nh}{c^2}\right)t}e^{-\left(\frac{nh}{c}\right)^2}f_s(nh)\sin(2\pi\nu t)dt
\\
&= -2i\frac{h}{\sqrt\pi}\int_{0}^\infty\sum_{n = -N}^N e^{-t^2 + \frac{2nh}{c}t}e^{-\left(\frac{nh}{c}\right)^2}f_s(nh)\sin(2\pi\nu ct)dt
\\
&= i\frac{h}{\sqrt\pi}\int_{0}^\infty\sum_{n = -N}^N e^{-\frac{t^2}{4} - \frac{nh}{c}t}e^{-\left(\frac{nh}{c}\right)^2}f_s(nh)\sin(\pi\nu ct)dt.
\end{aligned}
\]
Comparing this equation with the imaginary Voigt function \eqref{IVF}, we can rearrange equation above as
\begin{equation}
\label{FT4SF} 
{\hat f}_s(\nu) \approx ih\sum_{n = -N}^N e^{-\left(\frac{nh}{c}\right)^2}f_s(nh)L\left(\pi\nu c,\frac{nh}{c}\right).
\end{equation}

The applied methodology for the FTs of the rectangular and sawtooth \linebreak functions \eqref{FT4RF} and \eqref{FT4SF} can be generalized to any even or odd functions as
\begin{equation}
\label{FT4EF} 
{\hat f}_{even}(\nu) \approx h\sum_{n = -N}^N e^{-\left(\frac{nh}{c}\right)^2}f_{even}(nh)K\left(\pi\nu c,\frac{nh}{c}\right)
\end{equation}
and
\begin{equation}
\label{FT4OF} 
{\hat f}_{odd}(\nu) \approx ih\sum_{n = -N}^N e^{-\left(\frac{nh}{c}\right)^2}f_{odd}(nh)L\left(\pi\nu c,\frac{nh}{c}\right).
\end{equation}

Due to symmetric properties of the even and odd functions (see Equations \eqref{I4EF} and \eqref{I4OF}), the~number of the summation terms in these approximations can be reduced by a factor of~two.

Let us rearrange Equation \eqref{FT4EF} in the following form
\[
\begin{aligned}
{\hat f}_{even}(\nu) \approx & \, h\left[f_{even}(0)K(\pi\nu c,0) + \sum_{n = 1}^N e^{-\left(\frac{nh}{c}\right)^2}f_{even}(nh)K\left(\pi\nu c,\frac{nh}{c}\right)\right.
\\
&\left. + \sum_{n = -N}^{-1} e^{-\left(\frac{nh}{c}\right)^2}f_{even}(nh)K\left(\pi\nu c,\frac{nh}{c} \right)\right].
\end{aligned}
\]

Change of summation index~as
\[
\sum_{n = -N}^{- 1} e^{-\left(\frac{nh}{c}\right)^2}f_{even}(nh)K\left(\pi\nu c,\frac{nh}{c}\right) = \sum_{n = 1}^N e^{-\left(\frac{-nh}{c}\right)^2}f_{even}(-nh)K\left(\pi\nu c,\frac{- nh}{c}\right) 
\]
leads to
\[
\begin{aligned}
{\hat f}_{even}(\nu) \approx &\, h\left[f_{even}(0)K(\pi\nu c,0) + \sum_{n = 1}^N e^{-\left(\frac{nh}{c}\right)^2}f_{even}(nh)K\left(\pi\nu c,\frac{nh}{c}\right)\right.  
\\
&\left. + \sum_{n = 1}^N e^{-\left(\frac{-nh}{c}\right)^2}f_{even}(-nh)K\left(\pi\nu c,\frac{-nh}{c} \right)\right].
\end{aligned}
\]

Since according to Equation \eqref{I4EF}
\[
e^{-\left(\frac{-nh}{c}\right)^2} = e^{-\left(\frac{nh}{c}\right)^2},
\]
\[
f_{even}(-nh) = f_{even}(nh),
\]
and since
\[
K(x,0)=e^{-x^2} \Rightarrow K(\pi\nu c,0)=e^{-(\pi\nu c)^2},
\]
we can write
\begin{equation}
\label{FT4LT1} 
{\hat f}_{even}(\nu) \approx h \left(f_{even}(0)e^{-(\pi\nu c)^2} + \sum_{n = 1}^N {e^{-\left(\frac{nh}{c}\right)^2}f_{even}(nh)V_k\left(\pi\nu c,\frac{nh}{c}\right)}\right),
\end{equation}
where
\begin{equation}
\label{PVK} 
V_k\left(\pi\nu c,\frac{nh}{c}\right)=K\left(\pi\nu c,\frac{nh}{c}\right) + K\left(\pi\nu c,\frac{-nh}{c}\right)
\end{equation}
Equation \eqref{FT4OF} can be expanded as
\begin{equation}
\label{FT4OFE} 
\begin{aligned}
{\hat f}_{odd}(\nu) \approx &\, ih\left[f_{odd}(0)L(\pi\nu c,0) + \sum_{n = 1}^N e^{-\left(\frac{nh}{c} \right)^2}f_{odd}(nh)L\left(\pi\nu c,\frac{nh}{c}\right)\right.
\\
&+ \left.\sum_{n = -N}^{-1} e^{-\left(\frac{nh}{c}\right)^2}f_{odd}(nh)L\left(\pi\nu c,\frac{nh}{c} \right)\right].
\end{aligned}
\end{equation}

Taking into consideration that
\[
\sum_{n = -N}^{-1} e^{-\left(\frac{nh}{c}\right)^2}f_{odd}(nh)L\left(\pi\nu c,\frac{nh}{c}\right) = \sum_{n = 1}^N e^{-\left(\frac{-nh}{c}\right)^2}f_{odd}(-nh)L\left(\pi\nu c,\frac{- nh}{c}\right) 
\]
such that according to \eqref{I4OF}
\[
f_{odd}(-nh) = -f_{odd}(nh)
\]
and since
\[
L(\pi\nu c,0)=0,
\]
we can recast Equation \eqref{FT4OFE} as given by
\[
\begin{aligned}
{\hat f}_{odd}(\nu) \approx &\, ih\left[\sum_{n = 1}^N {e^{-\left(\frac{nh}{c}\right)^2}f_{odd}(nh)L\left(\pi\nu c,\frac{nh}{c}\right)}\right.
\\
&\left. + \sum_{n = 1}^N {e^{-\left(\frac{-nh}{c} \right)^2}f_{odd}(-nh)L\left(\pi\nu c,\frac{-nh}{c} \right)} \right]
\end{aligned}
\]
or
\[
\begin{aligned}
{\hat f}_{odd}(\nu) \approx &\, ih\Bigg[f_{odd}(0)L(\pi\nu c,0)
\\
&+ \sum_{n = 1}^N e^{-\left(\frac{nh}{c}\right)^2}f_{odd}(nh)\left(L\left(\pi\nu c,\frac{nh}{c}\right) - L\left(\pi\nu c,\frac{-nh}{c}\right)\right)\Bigg]
\end{aligned}
\]
or
\begin{equation}
\label{FT4LT2} 
{\hat f}_{odd}(\nu) \approx \, ih\sum_{n = 1}^N e^{-\left(\frac{nh}{c}\right)^2}f_{odd}(nh)V_\ell\left(\pi\nu c,\frac{nh}{c}\right),
\end{equation}
where
\begin{equation}
\label{PVL} 
V_\ell\left(\pi\nu c,\frac{nh}{c}\right)=L\left(\pi\nu c,\frac{nh}{c}\right) - L\left(\pi\nu c,\frac{-nh}{c}\right).
\end{equation}

At first glance, the FT Formulas \eqref{FT4LT1} and \eqref{FT4LT2} may appear computationally costly as they require user defined (external) function files. However, unlike Equation \eqref{RA}, the~expansion coefficients $V_k(\pi\nu c,nh/c)$ and $V_\ell(\pi\nu c,nh/c)$, shown by  Equations \eqref{PVK} and \eqref{PVL}, respectively, do not require re-computations every time when we change the sampled function $f(t)$ to any other function. Since the values $V_k(\pi\nu c,nh/c)$ and $V_\ell(\pi\nu c,nh/c)$ 
always remain the same regardless of the shape of the sampled function, these values can be precomputed in form of the look-up tables. Such implementation makes computation rapid as the required values $V_k(\pi\nu c,nh/c)$ and $V_\ell(\pi\nu c,nh/c)$ can be instantly picked up from the computer memory during computation of the FT. Furthermore, many algorithms, based on rational approximations, have been developed for rapid and high-accuracy computation of the Voigt/complex error function~\cite{Humlicek1979, Weideman1994, Wells1999, Letchworth2007, Schreier2011, Schreier2018, Abrarov2018a, Abrarov2018b, Azah2021, Thompson2024, Mangaldan}. Therefore, the~proposed technique can also be used as an alternative for rational approximation \eqref{RA} of the~FT.

\subsection{Numerical results}

The FT of the rectangular function can be readily found by integration as
\begin{equation}
\label{AFT1} 
{\hat f}_r(\nu) = \int_{-\infty}^\infty {f_r(t)e^{-2\pi i\nu t}dt = \int_{-1/2}^{1/2} e^{-2\pi i\nu t}dt = \, {\rm{sinc}}(\pi \nu)}.
\end{equation}
Similarly, we can find the FT of the sawtooth function as follows
\begin{equation}
\label{AFT2} 
\begin{aligned}
{\hat f}_s(\nu) &= \int_{-\infty}^\infty f_s(t)e^{- 2\pi i\nu t}dt = \int_{-\infty}^\infty tf_r(t)e^{-2\pi i\nu t}dt
\\
&= \int_{-1/2}^{1/2} te^{-2\pi i\nu t}dt = \, i\frac{\pi\nu\cos(\pi\nu) - \sin(\pi\nu)}{2\pi^2\nu^2},
\end{aligned}
\end{equation}
We can use these analytical results for comparison with numerical FTs computed by using Equations \eqref{FT4LT1} and \eqref{FT4LT2}.

Figure~3 shows numerical FTs of rectangular and the sawtooth functions by light green and magenta curves, respectively. Equations \eqref{AFT1} and \eqref{AFT2} are also shown by dashed black~curves.

\begin{figure}[ht]
\begin{center}
\includegraphics[width=24pc]{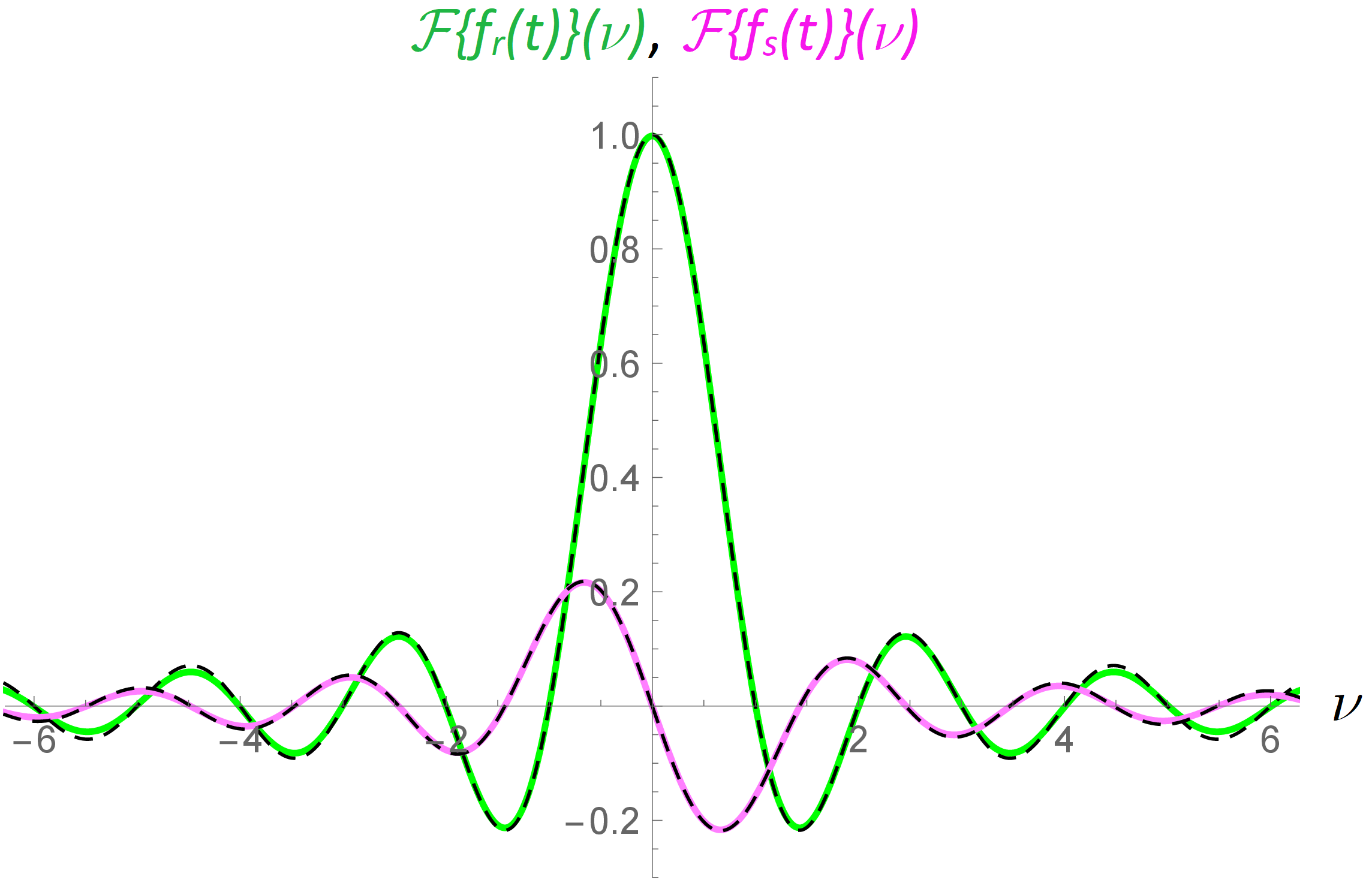}\hspace{2pc}%
\begin{minipage}[b]{28pc}
\vspace{0.3cm}
{\sffamily {\bf{Fig. 3.}} Fourier transforms of the rectangular and sawtooth functions (dashed black curves) and their approximations (green and magenta curves, respectively).}
\end{minipage}
\end{center}
\end{figure}

Figure~4 illustrates the absolute differences $\Delta_r$ and $\Delta_s$ between Equations \eqref{AFT1} and \eqref{AFT2} and numerical FT approximations \eqref{FT4LT1} and \eqref{FT4LT2} at $h=0.02$, $c=0.025$ and $N=25$ by blue and red curves, respectively. As~we can see from this numerical examples, despite the abrupt behavior of the rectangular and sawtooth functions $f_r(t)$ and $f_s(t)$ at $t_{1,2} = \pm 1/2$, the~FT approximations \eqref{FT4RF} and \eqref{FT4SF} can provide reasonable~accuracies.

\begin{figure}[ht]
\begin{center}
\includegraphics[width=24pc]{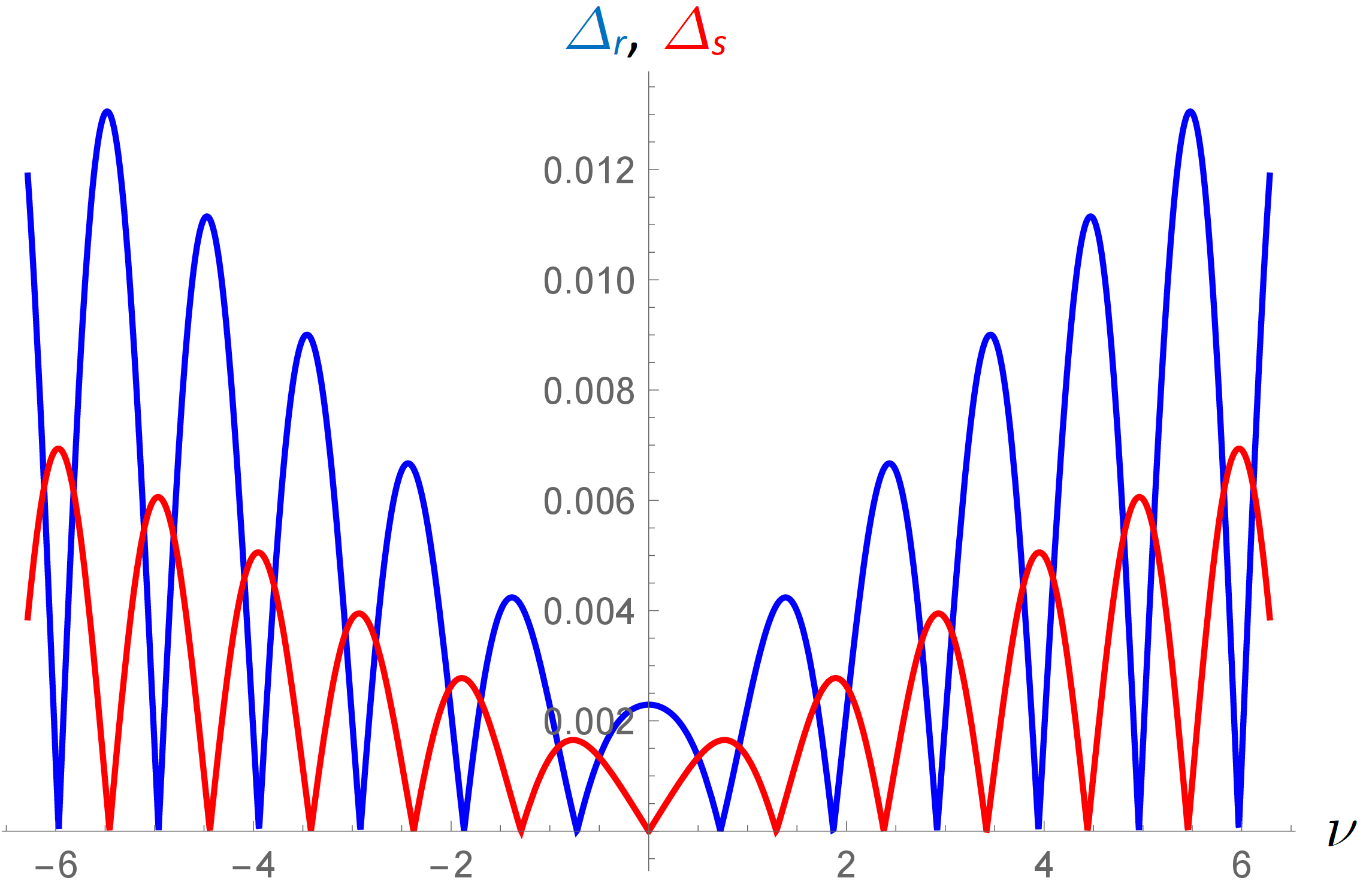}\hspace{2pc}%
\begin{minipage}[b]{28pc}
\vspace{0.3cm}
{\sffamily {\bf{Fig. 4.}} Absolute differences between equations \eqref{AFT1}, \eqref{AFT2} and their approximations. Blue curve corresponds to equation \eqref{AFT1}, red curve corresponds to equation \eqref{AFT2}.}
\end{minipage}
\end{center}
\end{figure}

The accuracy of the FT can be significantly better for the well-behaved functions. As~an example, consider a function
\begin{equation}
\label{WBF} 
g(t) = e^{-(6\pi t)^2} - \sin(32t)e^{-\left(7\pi t\right)^2}.
\end{equation}
The first term of this function is even
\begin{equation}
\label{ECGF} 
g_{even}(t) = \frac{g(t) + g(-t)}{2} = e^{-(6\pi t)^2}
\end{equation}
while its second term is odd
\begin{equation}
\label{OCGF} 
g_{odd}(t) = \frac{g(t) - g(-t)}{2} = -\sin(32t)e^{-(7\pi t)^2}.
\end{equation}

Figure 5 shows the functions $g(t)$, $g_{even}(t)$ and $g_{odd}(t)$ by blue, red, and green colors, respectively.

\begin{figure}[ht]
\begin{center}
\includegraphics[width=24pc]{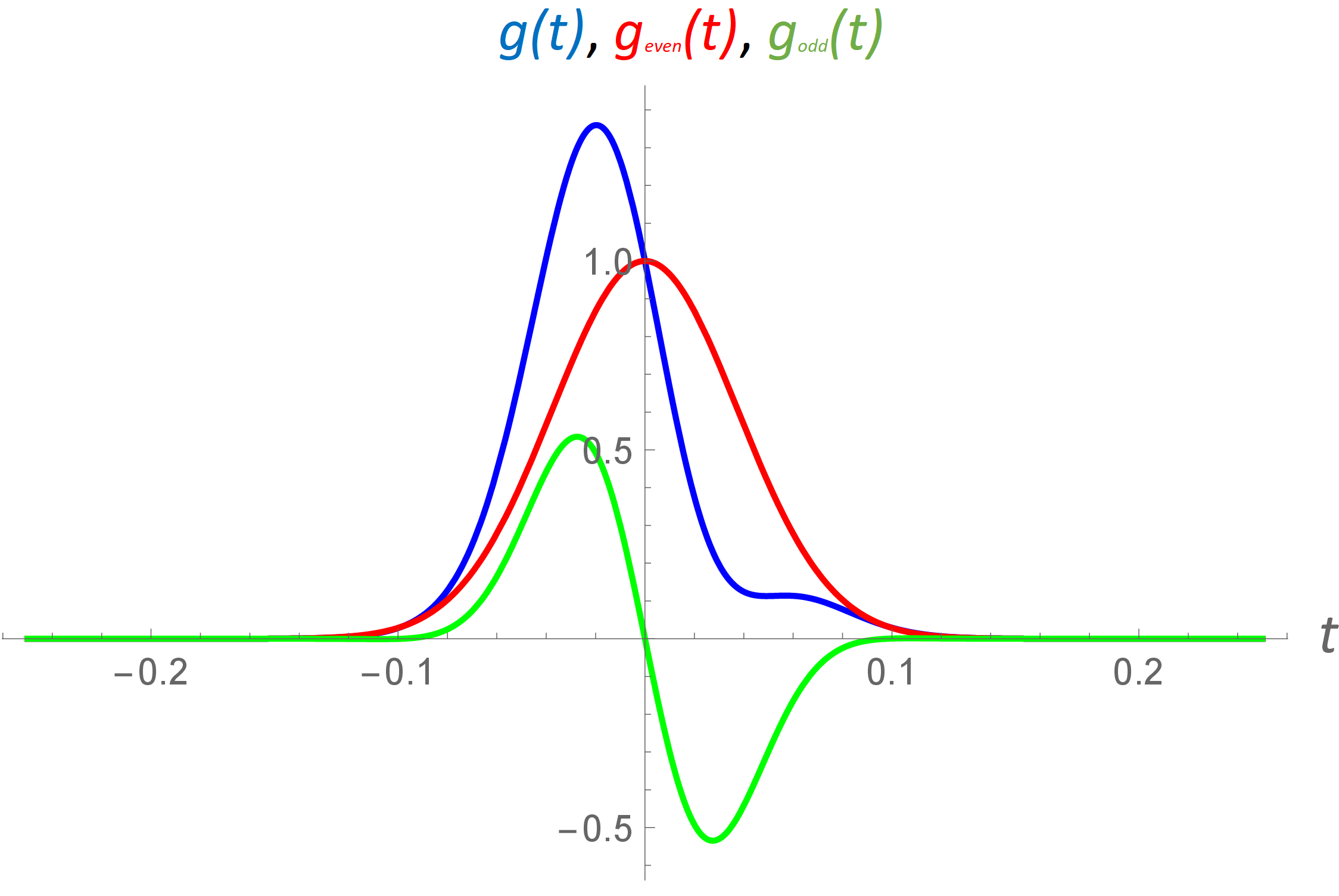}\hspace{2pc}%
\begin{minipage}[b]{28pc}
\vspace{0.3cm}
{\sffamily {\bf{Fig. 5.}} Equation \eqref{WBF} and its even and odd components (blue, red and green curves, respectively).}
\end{minipage}
\end{center}
\end{figure}

The FTs of even and odd components \eqref{ECGF} and \eqref{OCGF} can be obtained analytically. In~particular, substituting Equations \eqref{ECGF} and \eqref{OCGF} into FT Formulas \eqref{FT2} and \eqref{FT3} yields
\begin{equation}
\label{FT4EC} 
{\hat g}_{even}(\nu) = \frac{e^{-\left(\frac{\nu}{6}\right)^2}}{6\sqrt\pi}
\end{equation}
and
\begin{equation}
\label{FT4OC} 
{\hat g}_{odd}(\nu) = i\frac{1}{14\sqrt\pi}e^{-\left(\frac{16 + \pi\nu}{7\pi}\right)^2}\left(e^\frac{64\nu}{49\pi} - 1\right).
\end{equation}

Figure~6 depicts the absolute errors $\Delta_{even}$ and $\Delta_{odd}$ between Equations \eqref{FT4EC}, \eqref{FT4OC} and corresponding FT approximations \eqref{FT4LT1} and \eqref{FT4LT2} at $h=0.004$, $c=0.0045$, and $N=30$ by the blue and red colors, respectively. As~we can see, the~absolute differences do not exceed $0.00035$ and $0.0005$ for the functions and ${\hat g}_{even}(\nu)$ and ${\hat g}_{odd}(\nu)$, respectively, over~the wide interval of $\nu$.

It should be noted that further decreases in the values $h$, $c$ and increases in the integer $N$ can significantly improve the accuracy. For~example, in~our previous work~\cite{Abrarov2021} we have shown numerical examples where the FT of well-behavior functions can provide absolute difference $\sim 10^{-10}$.

\begin{figure}[ht]
\begin{center}
\includegraphics[width=24pc]{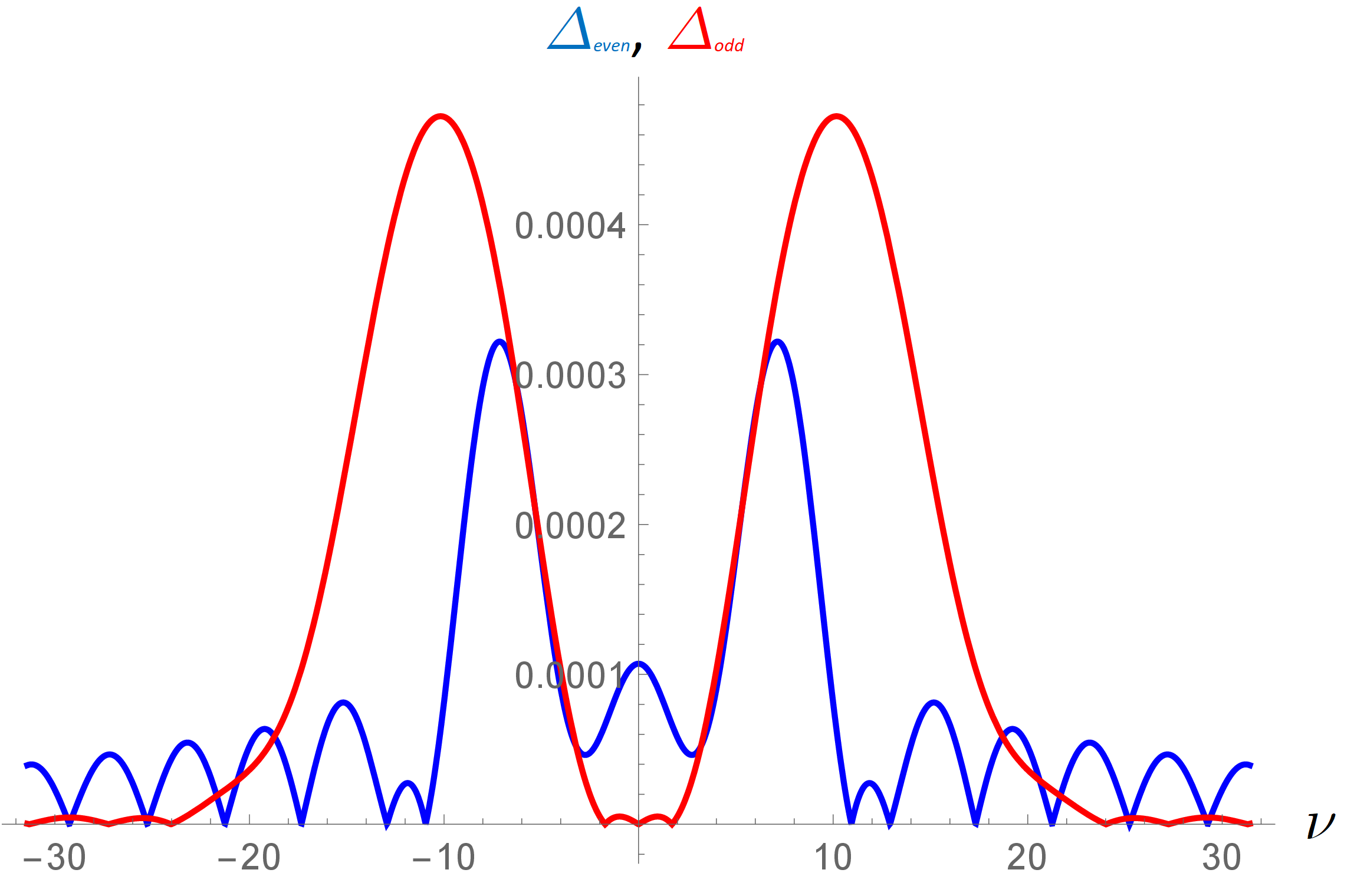}\hspace{2pc}%
\begin{minipage}[b]{28pc}
\vspace{0.3cm}
{\sffamily {\bf{Fig. 6.}} Absolute differences between equations \eqref{FT4LT1}, \eqref{FT4LT2} and their approximations. Blue curve corresponds to equation \eqref{FT4LT1}, red curve corresponds to equation \eqref{FT4LT2}.}
\end{minipage}
\end{center}
\end{figure}

\subsection{Trigonometric forms}

As previously mentioned,~Equations \eqref{FT4EF}, \eqref{FT4OF} and their variations \eqref{FT4LT1}, \eqref{FT4LT2} can be used to implement a rational approximation of the FT. We can also show how to represent these equations in trigonometric~forms.

The real part of Equation \eqref{I4CEF} gives
\[
K(x,-y) = 2e^{-x^2 + y^2}\cos(2xy) - K(x,y).
\]
Therefore, Equation \eqref{PVK} can be further simplified as
\[
\begin{aligned}
V_k\left(\pi\nu c,\frac{nh}{c}\right) &= K\left(\pi\nu c,\frac{nh}{c}\right) + 2e^{-(\pi\nu c)^2+\frac{nh}{c}}\cos\left(\pi\nu c,\frac{nh}{c}\right) - K\left(\pi\nu c,\frac{nh}{c}\right)
\\
&=2e^{-(\pi\nu c)^2+\frac{nh}{c}}\cos\left(\pi\nu c,\frac{nh}{c}\right)
\end{aligned}
\]
leading to
\begin{equation}
\label{FT4TF1} 
{\hat f}_{even}(\nu) \approx he^{-\left(\pi\nu c\right)^2}\left(f_{even}(0) + 2\sum_{n = 1}^N f_{even}(nh)\cos(2\pi\nu nh)\right)
\end{equation}
in accordance with Equation \eqref{FT4LT1}.

The imaginary part of Equation \eqref{I4CEF} provides
\[
L(x,-y) = 2e^{-x^2 + y^2}\sin(2xy) + L(x,y).
\]
Substituting this equation into approximation \eqref{PVL} yields
\[
\begin{aligned}
V_\ell\left(\pi\nu c,\frac{nh}{c}\right) &= L\left(\pi\nu c,\frac{nh}{c}\right) - 2e^{-(\pi\nu c)^2 + \left(\frac{nh}{c}\right)^2}\sin(2\pi\nu nh) - L\left(\pi\nu c,\frac{nh}{c}\right)
\\
&= -2e^{-(\pi\nu c)^2 + \left(\frac{nh}{c}\right)^2}\sin(2\pi\nu nh).
\end{aligned}
\]
This results in
\begin{equation}
\label{FT4TF2} 
{\hat f}_{odd}(\nu) \approx -2ihe^{-(\pi\nu c)^2}\sum_{n = 1}^N f_{odd}(nh)\left(\sin(2\pi\nu nh)\right)
\end{equation}
according to Equation \eqref{FT4LT2}.

As we can see, Equations \eqref{FT4TF1} and \eqref{FT4TF2} represent trigonometric versions of Equations \eqref{FT4LT1} and \eqref{FT4LT2}, respectively. Although~the objective of this work is a development of the FT based on rational approximation of the Voigt functions rather than an accelerated computation, it would be interesting to note that despite trigonometric representations, Equations \eqref{FT4TF1} and \eqref{FT4TF2} do not have any advantage over Equations \eqref{FT4LT1} and \eqref{FT4LT2} since all required values $V_k(\pi\nu c,nh/c)$ and $V_\ell(\pi\nu c,nh/c)$ as well as $\cos(2\pi\nu nh)$ and $\sin(2\pi\nu nh)$ can be precomputed and saved in computer memory in form of look-up tables. Therefore, both pairs of Equations \eqref{FT4LT1} and \eqref{FT4LT2} or \eqref{FT4TF1} and \eqref{FT4TF2} can be equally used with same computational speed and accuracy at given fitting parameters $h$, $c$, and $N$.

Combining Equations \eqref{EFI}, \eqref{OFI}, \eqref{FTI}, \eqref{FT4TF1}, and \eqref{FT4TF2} together, we obtain
\[
\begin{aligned}
{\hat f}(\nu) \approx &\, he^{-(\pi\nu c)^2}\Bigg[f(0) + \sum_{n = 1}^N \left(\left(f(nh) + f(-nh)\right)\cos(2\pi\nu nh)\right.
\\
&\left. -i\left(f(nh) - f(-nh)\right)\left(\sin(2\pi\nu nh)\right)\right)\Bigg]
\end{aligned}
\]
or
\begin{equation}
\label{TFTA} 
{\hat f}(\nu) \approx he^{-(\pi\nu c)^2}\left(f(0) + \sum_{n = 1}^N {\left(f(nh)e^{-2\pi i\nu nh} + f(-nh)e^{2\pi i\nu nh}\right)}\right).
\end{equation}

We can see that Equation \eqref{TFTA} resembles the discrete Fourier transfer (DFT) \cite{Hansen2014, Bracewell2000}. However, unlike the DFT, Equation \eqref{TFTA} provides non-periodic output (solitary pulse) due to exponential multiplier $e^{-(\pi\nu c)^2}$ acting like the Hamming window that sometimes may be introduced to eliminate periodicity in time or frequency domains. Remarkably, this exponential multiplier in Equation \eqref{TFTA} was not introduced but appeared naturally in~derivation.

Similar to the traditional DFT, the~real and imaginary parts of the functions $e^{-2\pi i\nu nh}$ and $e^{2\pi i\nu nh}$ in Equation \eqref{TFTA} are harmonic. Consequently, due to presence of these harmonics, we can construct the butterfly diagram for implementation of the fast Fourier transform (FFT) (see chapters 3--5 in~\cite{Hansen2014} showing how the functions $e^{-2\pi i\nu nh}$ and $e^{2\pi i\nu nh}$ can be used to implement FFT if we choose $1/h = N + 1$).

In Formula \eqref{TFTA}, the exponential multiplier $e^{-(\pi\nu c)^2}$ does not attenuate the higher frequency components in the FT when the values $h$, $c$, and $N$ are properly chosen (this keeps the shape of the solitary pulse undistorted). As~we can see from Figure~7a, with~$h = 0.004$ and $N = 30$ at $c = 0.001$ the absolute difference $\Delta_{even}$ does not exceed $0.00011$ (blue curve) while at $c = 0.002$ its value increases and exceeds $0.00014$ (red curve). However, this does not mean that the smaller value of the fitting parameter $c$ leads to a better result. The~higher frequency  components contribute to additional spikes in the function replica, calculated according to Equation \eqref{TFTA}. These equidistantly spanned spikes can be seen from Figure~7b illustrating the function replica at $c = 0.001$ and $c = 0.002$ by blue and red curves, respectively. Generally, the~best accuracy without undesirable spikes may be obtained when both values $c$ and $h$ are small enough and close to each other and when the integer $N$ is sufficiently~large.

\begin{figure}[ht]
\begin{center}
\includegraphics[width=32pc]{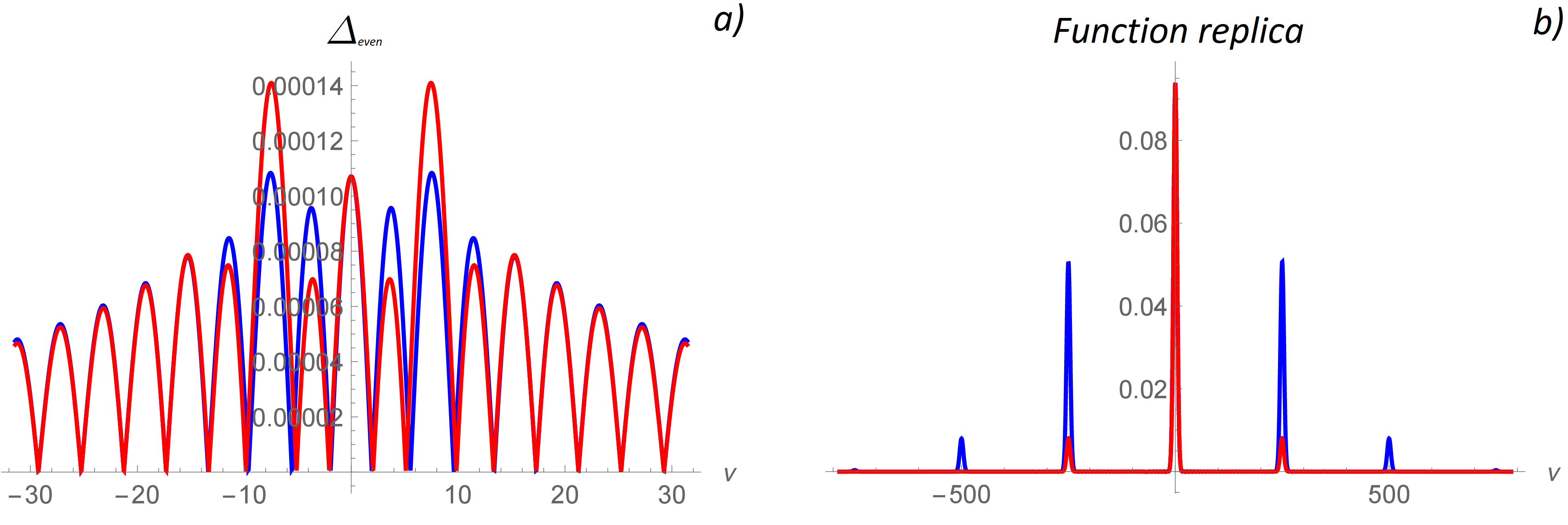}\hspace{2pc}%
\begin{minipage}[b]{28pc}
\vspace{0.3cm}
{\sffamily {\bf{Fig. 7.}} Absolute differences {\it a)} and function replica {\it b)} at $c = 0.001$ (blue curves) and $c = 0.002$ (red curves). The parameters $h = 0.004$ and $N = 30$ are the same.}
\end{minipage}
\end{center}
\end{figure}

Accuracy analysis can be focused on optimization of the fitting parameters $h$, $c$ and $N$. Our empirical results show that the Gibbs phenomenon~\cite{Hewitt1979, Jerri1998}, occurring as a result of sampling (see oscillations in Figure~2
), decreases when $h \to c$. However, this does not imply that the best optimization in the rational approximation of the FT can be achieved precisely at $h = c$ since, apart from the sampling, this process also involves integration with respect to the variable $\nu$. Therefore, absolute differences $\Delta_{even}$, $\Delta_{odd}$ for the rational  approximation of the FT are functions of four parameters $h$, $c$, $N$, and $\nu$. A~slight difference between two small parameters $h$ and $c$ gives us some flexibility in optimization. Furthermore, more complete accuracy analysis should also include consideration of the rational approximations that are used for computing the Voigt functions. Highly accurate algorithms for computation of the Voigt functions require a combination of several different methods of approximations for efficient domain-to-domain coverage over the complex plane (see~\cite{Humlicek1979, Weideman1994, Wells1999, Letchworth2007, Schreier2011, Schreier2018, Abrarov2018a, Abrarov2018b, Azah2021, Thompson2024} and literature therein for implementation of different methods of rational approximations). Consequently, optimization and more rigorous numerical analysis on accuracy requires a separate consideration and, therefore, it 
is beyond the scope of this~work.

Thus, based on the results described above the following key results obtained in this study may be highlighted~as
\begin{itemize}
\item It is shown that the FT can be expressed in terms of the Voigt functions;
\item Expansion coefficients in rational approximation of the FT always remains the same;
\item	Precomputed values of the Voigt functions can be stored in a computer memory in form of the look-up tables;
\item Application of the precomputed values in the look-up tables can be used for more efficient computation of the FT;
\item There is a trigonometric form of the FT based on the Voigt functions.
\end{itemize}

Figure~8 illustrates a flowchart of the FT by using the look-up tables that can be computed with help of the Voigt functions. Initially, an~input function $f(t)$ undergoes sampling. After~that the discrete sampled function is passed to the module in which the FT occurs. The~FT module is linked to the look-up tables and the entire process can be controlled by providing the corresponding fitting parameters. Once the FT is completed, its signal is still discrete. Therefore, at~the final stage the interpolation (if required) can be used to fill the gaps between all adjacent points. After~interpolation the output function becomes continues. The~abbreviations CF and DF shown in the  flowchart are for the continuous and discrete functions, respectively.

\begin{figure}[ht]
\begin{center}
\includegraphics[width=32pc]{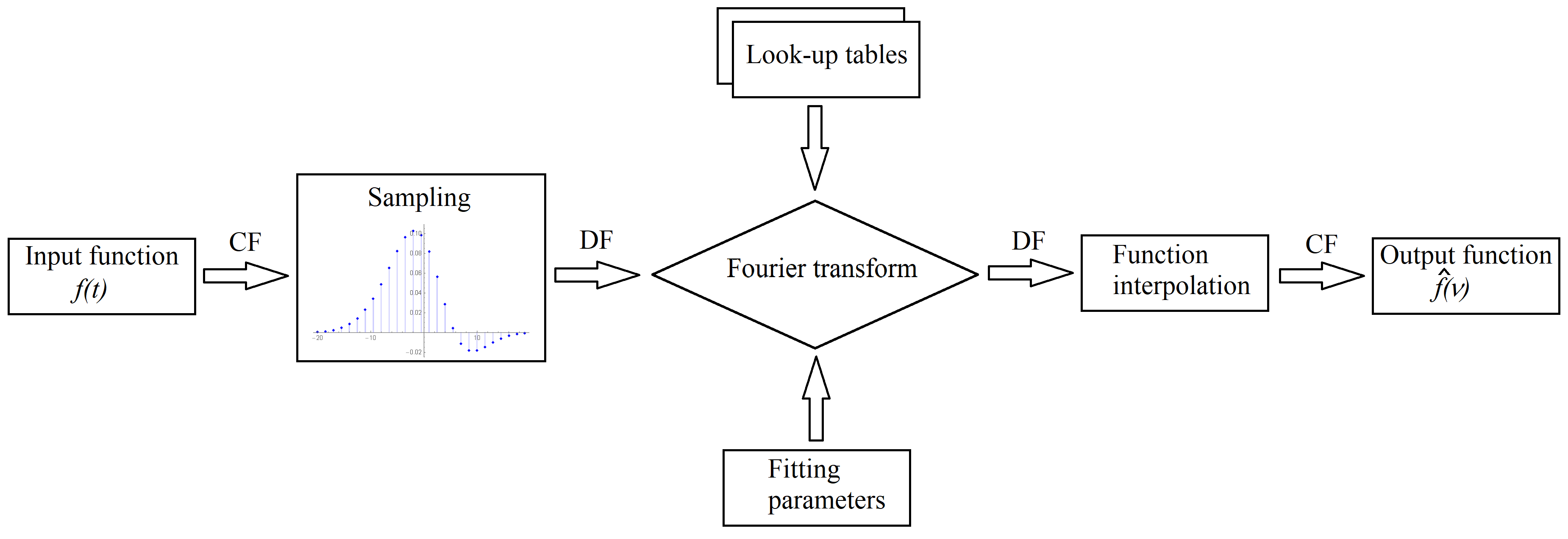}\hspace{2pc}%
\begin{minipage}[b]{28pc}
\vspace{0.3cm}
{\sffamily {\bf{Fig. 8.}} Flowchart of the FT by using the look-up tables that can be precomputed with help of the Voigt functions. Input and output are continuous functions (CFs). Discrete functions (DFs) are used in the intermediate steps of the~FT.}
\end{minipage}
\end{center}
\end{figure}

Examples of the numerical FTs based on equations \eqref{FT4LT1} and \eqref{FT4LT2} can be validated by running the Mathematica codes, shown in the next section.

\pagestyle{empty}
\section{Mathematica codes}

\subsection{Complex error function}

To compute the complex error function $w(z)$, we can use, for example, the algorithm provided in our work \cite{Abrarov2018a}. To cover the entire complex plane with high accuracy of computation, we used only three approximations. We will use four cells to define each equation and to distribute them accordingly on the complex plane.

The first cell below is to define equation (7) from \cite{Abrarov2018a}.
\small
\begin{verbatim}
Clear[a0, b0, c0, \[CapitalOmega], w1];

(* Fitting parameters *)

H = 0.25; \[Stigma] = 2.75; M = 25; maxN = 23;

(* Expansion coefficients, 1st set *)

a0[m_] := a0[m] = ((Sqrt[Pi]*(m - 1/2))/(2*M^2*H))*
  Sum[E^(\[Stigma]^2/4 - n^2*H^2)*Sin[(Pi*(m - 1/2)*
    (n*H + \[Stigma]/2))/(M*H)], {n, -maxN, maxN}];

b0[m_] := b0[m] = (-(I/(M*Sqrt[Pi])))*
  Sum[E^(\[Stigma]^2/4 - n^2*H^2)*Cos[(Pi*(m - 1/2)*
    (n*H + \[Stigma]/2))/(M*H)], {n, -maxN, maxN}];

c0[m_] := c0[m] = (Pi*(m - 1/2))/(2*M*H);

(* Equation (7) from Ref. [29] *)  

\[CapitalOmega][z_] := \[CapitalOmega][z] = Sum[(a0[m] + b0[m]*z)/
  (c0[m]^2 - z^2), {m, 1, M - 2}];

w1[z_] := \[CapitalOmega][z + I*(\[Stigma]/2)];
\end{verbatim}
\normalsize

The second cell is required to instantiate equation (8) from \cite{Abrarov2018a}.
\small
\begin{verbatim}
Clear[a1, b1, c1, d1, w2];

(* Expansion coefficients, 2nd set *)

a1[m_] := a1[m] = b0[m]*(((Pi*(m - 1/2))/(2*M*H))^2 -
  (\[Stigma]/2)^2) + I*a0[m]*\[Stigma];

b1[m_] := b1[m] = b0[m];

c1[m_] := c1[m] = (((Pi*(m - 1/2))/(2*M*H))^2 +
  (\[Stigma]/2)^2)^2;

d1[m_] := d1[m] = 2*((Pi*(m - 1/2))/(2*M*H))^2 -
  \[Stigma]^2/2;

(* Equation (8) from Ref. [29] *)

w2[z_] := E^(-z^2) + z*Sum[(a1[m] - b1[m]*z^2)/
  (c1[m] - d1[m]*z^2 + z^4), {m, 1, M - 2}];
\end{verbatim}
\normalsize

This third cell is required to define equation (9) from \cite{Abrarov2018a}.
\small
\begin{verbatim}
Clear[w3];

(* Equation (9) from Ref. [29] *)
w3[z_] := I/(Sqrt[Pi]*(z - 1/(2*(z - 1/(z - 3/
  (2*(z - 2/(z - 5/(2*(z - 3/(z - 7/(2*(z - 4/
    (z - 9/(2*(z - 5/(z - 11/(2*z))))))))))))))))));
\end{verbatim}
\normalsize

Once these three equations are instantiated, we need to distribute them correspondingly on the complex plane. The forth cell below contains code to accomplish this task.
\small
\begin{verbatim}
Clear[wUp, w];

(* Complex error function for upper complex plane *)

wUp[z_] := If[Abs[z] > 8, w3[z],
  If[Im[z] > 0.05*Abs[Re[z]], w1[z], w2[z]]];

(* Complex error function for entire complex plane *)

w[z_] := If[Im[z] >= 0, wUp[z],
  Conjugate[2*E^-Conjugate[z]^2 - wUp[Conjugate[z]]]];
\end{verbatim}
\normalsize
Now the code for computation of the complex error function is ready to use.

It should be noted that another Mathematica code for high-accuracy computation of the complex error function can be downloaded from \cite{Mangaldan}. This code was written by Jan Mangaldan on the bases of three rational approximations described in our publication \cite{Abrarov2018b}.

\subsection{Fourier transform}

The Mathematica codes below consist of six cells. The code shown in the first cell below defines the Voigt function $K(x,y)$ and imaginary Voigt function $L(x,y)$, the rectangular function $f_r(t)$ and the sawtooth function $f_s(t)$.
\small
\begin{verbatim}
Clear[K,L,fr,fs];

(* Defining K(x,y) and L(x,y) functions *)

K[x_, y_] := Re[w[x + I*y]];
L[x_, y_] := Im[w[x + I*y]];

(* Rectangular function *)

fr[t_] := 1/((2*t)^(2*35) + 1);

(* Sawtooth function *)

fs[t_] := t*fr[t];
\end{verbatim}
\normalsize

The second cell below generates two look-up tables for the values $V_k(\pi\nu c,nh/c)$ and $V_\ell(\pi\nu c,nh/c)$. It also generates a list of grid-points for the parameter $\nu$.
\small
\begin{verbatim}
Clear[lookUpTab1, lookUpTab2, nuList];

(* Parameters for computation *)

h = 0.02; c = 0.025; nMax = 25;

(* Computing two look-up tables *)

lookUpTab1 = Table[{K[Pi*\[Nu]*c, (n*h)/c] +
  K[Pi*\[Nu]*c, ((-n)*h)/c]}, {n, 1, nMax},
    {\[Nu], -2*Pi, 2*Pi, 0.1}];

lookUpTab2 = Table[{L[Pi*\[Nu]*c, (n*h)/c] -
  L[Pi*\[Nu]*c, ((-n)*h)/c]}, {n, 1, nMax},
    {\[Nu], -2*Pi, 2*Pi, 0.1}];

nuList = Table[\[Nu], {\[Nu], -2*Pi, 2*Pi, 0.1}];
\end{verbatim}
\normalsize

The code in the following third cell is needed to format two look-up table with values $V_k(\pi\nu c,nh/c)$ and $V_\ell(\pi\nu c,nh/c)$ for plotting the graphs. It is also required to join $V_k(\pi\nu c,nh/c)$ and $V_\ell(\pi\nu c,nh/c)$ values together with corresponding values of the parameter $\nu$.

\small
\begin{verbatim}
Clear[ftList1, ftList2];

(* Main computations by using look up tables *)

ftList1 = Flatten[h*(fr[0]/E^(Pi*nuList*c)^2 +
  Sum[(fr[n*h]*lookUpTab1[[n]])/E^((n*h)/c)^2, {n, 1, nMax}])];

ftList2 = Flatten[h*Sum[(fs[n*h]*lookUpTab2[[n]])/E^((n*h)/c)^2,
  {n, 1, nMax}]];

(* Arranging FT data lists *)

ftList1 = Table[{nuList[[n]], ftList1[[n]]}, {n, 1, Length[nuList]}];

ftList2 = Table[{nuList[[n]], ftList2[[n]]}, {n, 1, Length[nuList]}];
\end{verbatim}
\normalsize

The code in next forth cell is required to generate the references in accordance with equations \eqref{AFT1} and \eqref{AFT2}.
\small
\begin{verbatim}
Clear[ftRef1, ftRef2];

(* FT references *)

ftRef1 = Table[{\[Nu], Sinc[Pi*\[Nu]]}, {\[Nu], -2*Pi, 2*Pi, 0.1}];

ftRef2 = Table[{\[Nu], (Pi*\[Nu]*Cos[Pi*\[Nu]] - Sin[Pi*\[Nu]])/(2*
  Pi^2*\[Nu]^2)}, {\[Nu], -2*Pi, 2*Pi, 0.1}];
\end{verbatim}
\normalsize

The code in the fifth cell below produces the graph shown in the Fig. 3.
\small
\begin{verbatim}
(* Plotting the graphs from the data lists *)

ListPlot[{ftList1, ftRef1, ftList2, ftRef2}, PlotRange -> All,
  Joined -> True, PlotStyle -> {{Lighter[Green, 0], Thickness[0.005]},
    {Black, Dashed, Thickness[0.0025]}, {Lighter[Magenta, 0.5],
Thickness[0.005]}, {Black, Dashed, Thickness[0.0025]}},
  PlotRange -> {{-2*Pi, 2*Pi}, {-0.3, 1.1}},
    AxesLabel -> {"\[Nu]", None}]
\end{verbatim}
\normalsize

Lastly, the code in following sixth cell applies the derived formula \eqref{TFTA} to generate the same Fig. 3 without Voigt functions.

\small
\begin{verbatim}
(* Plotting the graphs by using the FT formula (47) *)

f[t_] := fr[t] + fs[t];

ft[\[Nu]_] := (h*(f[0] + Sum[f[n*h]/E^(2*Pi*I*\[Nu]*n*h) +
  f[(-n)*h]*E^(2*Pi*I*\[Nu]*n*h), {n, 1, nMax}]))/
    E^(Pi*\[Nu]*c)^2;

Plot[{Re[ft[\[Nu]]], Sinc[Pi*\[Nu]], Im[ft[\[Nu]]],
  (Pi*\[Nu]*Cos[Pi*\[Nu]] - Sin[Pi*\[Nu]])/(2*(Pi*\[Nu])^2)},
    {\[Nu], -2*Pi, 2*Pi}, PlotRange -> All, PlotStyle ->
      {{Lighter[Green, 0], Thickness[0.005]}, {Black, Dashed,
Thickness[0.0025]}, {Lighter[Magenta, 0.5],
  Thickness[0.005]},{Black, Dashed, Thickness[0.0025]}},
    PlotRange -> {{-2*Pi, 2*Pi}, {-0.3, 1.1}},
      AxesLabel -> {"\[Nu]", None}]
\end{verbatim}
\normalsize

The codes shown in this section can be copy-pasted directly to the Mathematica notebook.

\pagestyle{plain}
\section{Application}

In this section, we show how to apply practically the results obtained above. Making change of the variable $f \to \hat{f}^{-1}$ in Equation \eqref{TFTA}, we have
\small
\begin{equation}
\label{FLSE1} 
f(t) \approx he^{-(\pi ct)^2}\left(\hat{f}^{-1}(0) + \sum_{n = 1}^N \left(\hat{f}^{-1}(nh)e^{-2\pi inht}+\hat{f}^{-1}(-nh)e^{2\pi inht}\right)\right).
\end{equation}
\normalsize

Approximation \eqref{FLSE1} requires complex numbers in computation. Therefore, its application is not optimal. However, this approximation can be significantly simplified. In~particular, according to Equations \eqref{FT4TF1} and \eqref{FT4TF2}, the~series expansion \eqref{FLSE1} can also be expressed in form
\footnotesize
\begin{equation}
\label{FLSE2} 
\begin{aligned}
f(t) \approx &
\\
&\, he^{-(\pi ct)^2}\left[\hat{f}^{-1}_{even}(0) + 2\sum_{n = 1}^N \left(\hat{f}^{-1}_{even}(nh)\cos(2\pi tnh) - i\hat{f}^{-1}_{odd}(nh)\sin(2\pi tnh)\right)\right].
\end{aligned}
\end{equation}
\normalsize

Since
$$
\hat{f}^{-1}_{even}(\nu) = \Re\left\{\hat{f}^{-1}(\nu)\right\}
$$
and
$$
\hat{f}^{-1}_{odd}(\nu) = i\Im\left\{\hat{f}^{-1}(\nu)\right\},
$$
from Equation \eqref{FLSE2} it follows that
\footnotesize
\begin{equation}
\label{FLSE3} 
\begin{aligned}
f&(t) \approx he^{-(\pi ct)^2}
\\
&\times\left[\Re\left\{\hat{f}^{-1}(0)\right\} + 2\sum_{n = 1}^N \left(\Re\left\{\hat{f}^{-1}(nh)\right\}\cos(2\pi tnh) + \Im\left\{\hat{f}^{-1}(nh)\right\}\sin(2\pi tnh)\right)\right].
\end{aligned}
\end{equation}
\normalsize

Consider as an example a real valued function
\[
u(t) = \Re\left\{g(t)\right\} - \Im\left\{g(t)\right\}
\]
that according to Equations \eqref{FT4EC} and \eqref{FT4OC} is given by
\begin{equation}
\label{D4FU} 
u(t) = \frac{e^{-\left(\frac{t}{6}\right)^2}}{6\sqrt\pi} - \frac{e^{-\left(\frac{\pi t + 16}{7\pi}\right)^2} \left(e^{\frac{64t}{49\pi}} - 1\right)}{14\sqrt\pi}.
\end{equation}
Consequently, the~IFT of the function \eqref{D4FU} can be found by substituting this function into Equation \eqref{IFT}. This substitution results in the following IFT
\begin{equation}
\label{FT4U} 
\hat{u}^{-1}(\nu) = e^{-(6\pi\nu)^2} - ie^{-(7\pi\nu)^2}\sin(32\nu).
\end{equation}

Finally, substituting Equation \eqref{FT4U} into series expansion \eqref{FLSE3}, we get 
\footnotesize
\begin{equation}
\label{E4SE} 
\begin{aligned}
u(t) &\approx he^{-(\pi ct)^2}
\\
&\times\left[1 + 2\sum_{n = 1}^N \left(e^{-(6\pi nh)^2}\cos(2\pi tnh) - e^{-(7\pi nh)^2}\sin(32 nh)\sin(2\pi tnh)\right)\right].
\end{aligned}
\end{equation}
\normalsize
\normalsize

Figure~9a shows Equation \eqref{D4FU} and its approximation \eqref{E4SE} by light blue and dashed-black curves respectively, computed at $h = 0.004$, $c = 0.004$ and $N = 30$. Figure~9b illustrates absolute difference $\Delta_{u(t)}$ between Equation \eqref{D4FU} and its approximation \eqref{E4SE}. From~these results, we can see that a reasonable accuracy can be obtained even at a relatively small numbers of the summation terms. Specifically, with~only $N = 30$ the absolute difference $\Delta_{u(t)}$ does not exceed $0.0006$.

\begin{figure}[ht]
\begin{center}
\includegraphics[width=32pc]{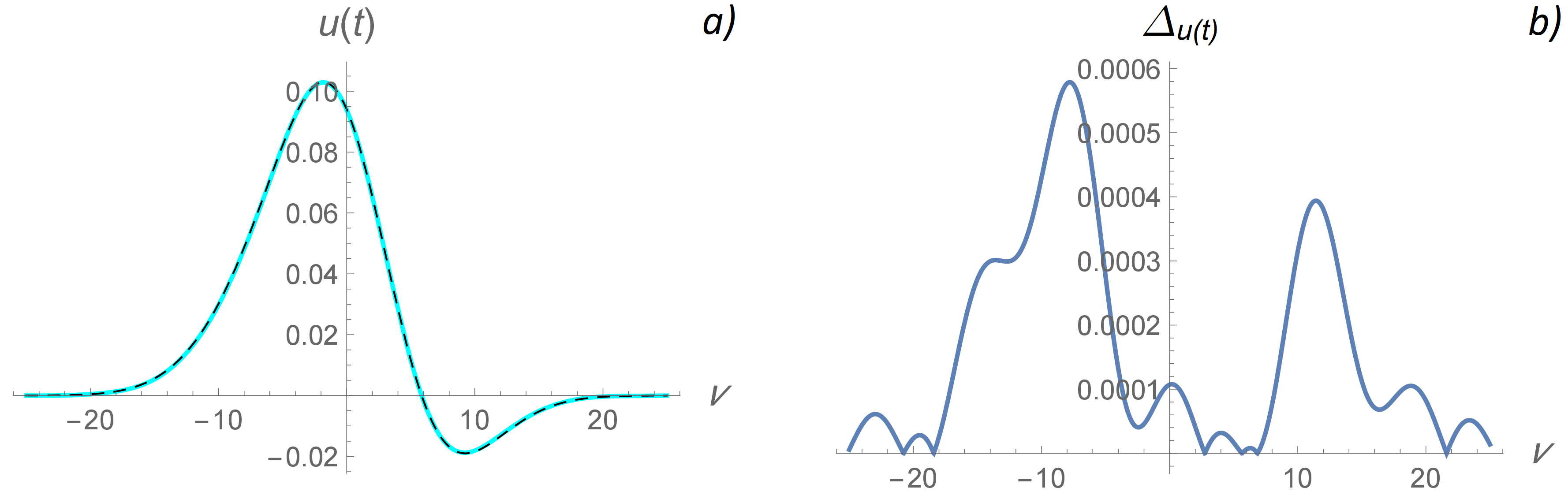}\hspace{2pc}%
\begin{minipage}[b]{28pc}
\vspace{0.3cm}
{\sffamily {\bf{Fig. 9.}} {\it a)} Equation \eqref{D4FU} and its approximation \eqref{E4SE}, shown by light blue and dashed-black curves, respectively. (\textbf{b}) Absolute difference $\Delta_{u(t)}$ between Equation \eqref{D4FU} and its approximation \eqref{E4SE}. The~following fitting parameters are used: $h = 0.004$, $c = 0.004$ and $N = 30$.}
\end{minipage}
\end{center}
\end{figure}

Equation \eqref{FLSE3} can be used as an efficient alternative to the traditional Fourier expansion series. In~contrast to the Fourier expansion series, this approach is more flexible as it can be used to expand a solitary pulse. Moreover, by~taking $c = 0$ we can use Equation \eqref{FLSE3} as a complete analog to the traditional Fourier expansion~series.

Apart from the absolute difference, we can also evaluate accuracy of \mbox{approximation \eqref{E4SE}} by using the relative error that can be defined as
\[
\varepsilon = \frac{{\rm Approximation} - {\rm Reference}}{\rm Reference}.
\]

\begin{figure}[ht]
\begin{center}
\includegraphics[width=18pc]{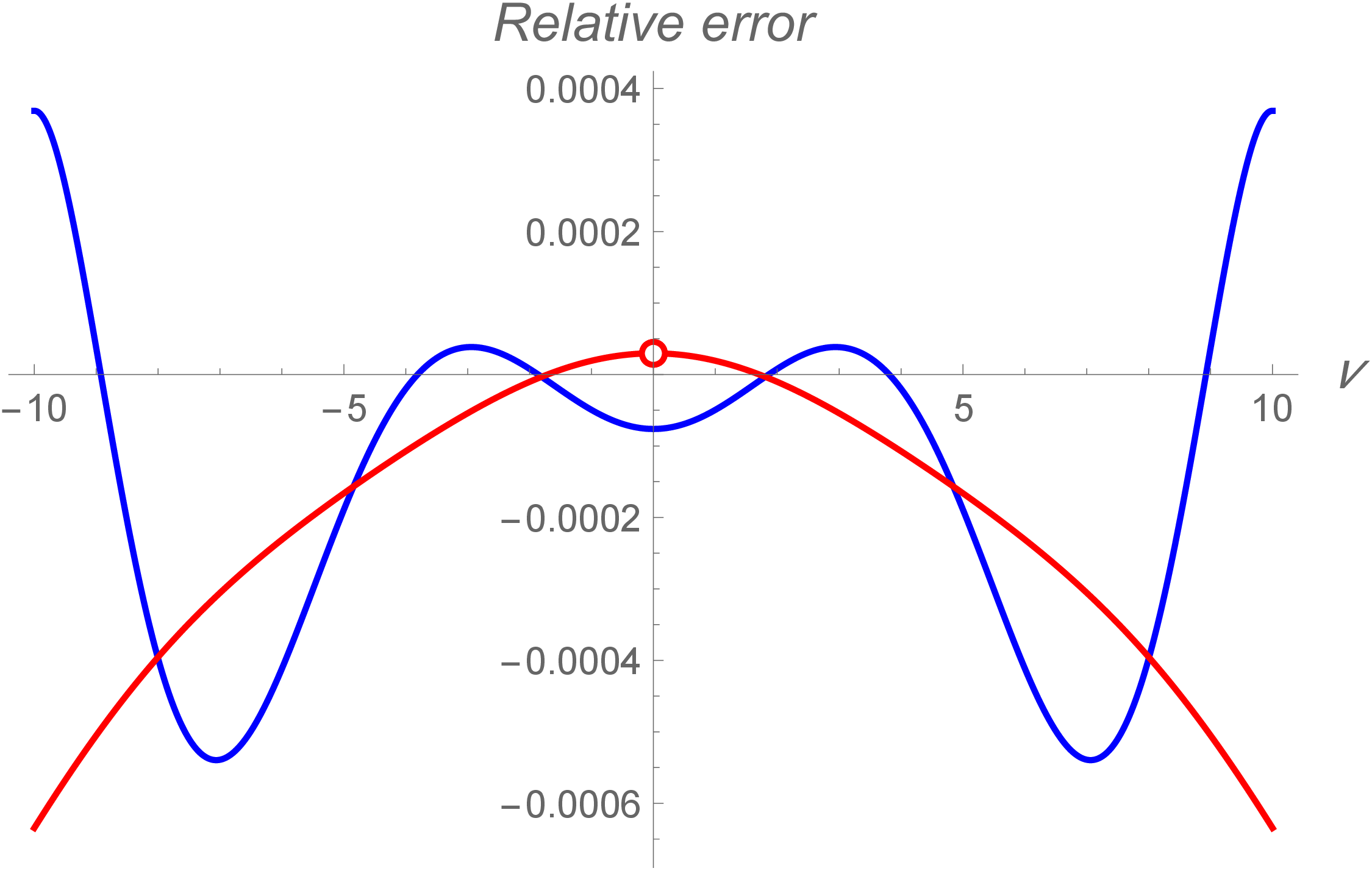}\hspace{2pc}%
\begin{minipage}[b]{28pc}
\vspace{0.3cm}
{\sffamily {\bf{Fig. 10.}} Relative error $\varepsilon$ for the even (blue) and odd (red) terms of Equation \eqref{D4FU} computed at $h = 0.0008, c = 0.0008$ and $N = 185$. Open circle located at center of the red curve shows the uncertainty point appearing due to division to~zero.}
\end{minipage}
\end{center}
\end{figure}

Figure~10 shows the relative error $\varepsilon$ for the even and odd terms of Equation \eqref{D4FU} by blue and red curves, respectively. As~we can see from this figure, within~the domain $t \in [-10,10]$ at $h = 0.0008, c = 0.0008$ and $N = 185$ the relative error $\varepsilon$ is smaller than the corresponding values of the original function \eqref{D4FU} by several orders of the magnitude. Although~beyond this domain the relative error rapidly increases, it is possible to extend the domain coverage by  further reducing the parameters $h$, $c$ and increasing the integer $N$. Therefore, it is a trade-off between the domain coverage and number of the summations terms $N$.

There is a point of uncertainty shown in Figure~10 by red empty circle located just above 
the origin. This point of uncertainty in $\varepsilon$ appears due to division to zero as the 
odd term of Equation \eqref{D4FU} equals zero when $t = 0$.

\section{Conclusion}

An alternative method of rational approximation of the FT based on the real and imaginary parts of the complex error function \eqref{CEF1} is developed. Unlike the rational approximation \eqref{RA} of the FT, the~expansion coefficients $V_k(\pi\nu c, nh/c)$ and $V_\ell(\pi\nu c, nh/c)$ in this method do not depend on values of the sampled function $f(t)$. Since the values of the Voigt functions always remain the same, this approach can be used for rapid computation with help of look-up tables. We also show that this rational approximation of the FT can also be rearranged in a trigonometric form \eqref{TFTA} with an exponential multiplier $e^{-(\pi\nu c)^2}$ acting like the Hamming window that removes~periodicity.

\section*{Acknowledgment}

This work was supported by National Research Council Canada, Thoth Technology Inc., York University and Epic College of Technology. The authors wish to thank anonymous reviewers for their constructive comments and~recommendations.


\end{document}